%% file: C.tex
\input JOSEPH

\global\ssecno=0\global\secno=0
\def\page#1#2{\leaders\hbox to 2mm{\hfil.\hfil}\hfill
\rlap{\hbox to 15mm{\hfill #1}{\ \  p.#2}}\par
}

\def\tpage #1{\leaders\hbox to 2mm{\hfil.\hfil}\hfill
\rlap{\hbox to 11mm{\hfill #1}}\par
}
\def\npage {\vfill\eject \global\notenumber=1}

\input C0.tex
\npage
\input C1.tex

\input Ref.tex
\npage
\input C2.tex
\npage

\input Mat.tex
\end

%% file: JOSEPH.tex
\input option_keys

\input SSIND
\input FONTESMATH

\input MATH

\input pdffix.tex

%% file: option_keys
\ifx\optionkeymacros\undefined\else \fi

\catcode`\Œ=\active\defŒ{{\aa}}       
\catcode`\º=\active\defº{\int}        
\catcode`\=\active\def{\c c}        
\catcode`\¶=\active\def¶{\partial}    
\catcode`\Ä=\active\defÄ{\oint}       
\catcode`\Æ=\active\defÆ{\triangle}   
\catcode`\Â=\active\defÂ{\neg}        
\catcode`\µ=\active\defµ{\mu}         
\catcode`\¿=\active\def¿{{\o}}        
\catcode`\¹=\active\def¹{\pi}         
\catcode`\Ï=\active\defÏ{{\oe}}       
\catcode`\§=\active\def§{{\ss}}       
\catcode`\ =\active\def {\dagger}     
\catcode`\Ã=\active\defÃ{\sqrt}       
\catcode`\·=\active\def·{\Sigma}      
\catcode`\Å=\active\defÅ{\approx}     
\catcode`\½=\active\def½{\Omega}      
\catcode`\£=\active\def£{{\it\$}}     
\catcode`\°=\active\def°{\infty}      
\catcode`\¤=\active\def¤{{\S}}        
\catcode`\¦=\active\def¦{{\P}}        
\catcode`\¥=\active\def¥{\bullet}     
\catcode`\»=\active\def»{\leavevmode\raise.585ex\hbox{\b a}}      
\catcode`\¼=\active\def¼{\leavevmode\raise.6ex\hbox{\b o}}        
\catcode`\­=\active\def­{\not=}       
\catcode`\²=\active\def²{\leq}        
\catcode`\³=\active\def³{\geq}        
\catcode`\Ö=\active\defÖ{\div}        
\catcode`\É=\active\defÉ{{\dots}}     
\catcode`\¾=\active\def¾{{\ae}}       
\catcode`\Ç=\active\defÇ{\ll}         
\catcode`\Ò=\active\defÒ{``}          
\catcode`\Á=\active\defÁ{!`}          
\catcode`\¢=\active\def¢{\rlap/c}     


\catcode`\=\active\def{{\AA}}       
\catcode`\'=\active\def'{\c C}        
\catcode`\¯=\active\def¯{{\O}}        
\catcode`\¸=\active\def¸{\Pi}         
\catcode`\Î=\active\defÎ{{\OE}}       
\catcode`\®=\active\def®{{\AE}}       
\catcode`\×=\active\def×{\diamond}    
\catcode`\¡=\active\def¡{\accent'27}  
\catcode`\Ó=\active\defÓ{''}          
\catcode`\±=\active\def±{\pm}         
\catcode`\È=\active\defÈ{\gg}         
\catcode`\À=\active\defÀ{?`}          
\catcode`\Ð=\active\defÐ{--}          
\catcode`\Ñ=\active\defÑ{---}         


\catcode`\Š=\active\defŠ{\"a}        
\catcode`\'=\active\def'{\"e}        
\catcode`\•=\active\def•{\"{\i}}     
\catcode`\š=\active\defš{\"o}        
\catcode`\Ÿ=\active\defŸ{\"u}        
\catcode`\Ø=\active\defØ{\"y}        
\catcode`\€=\active\def€{\"A}        
\catcode`\…=\active\def…{\"O}        
\catcode`\†=\active\def†{\"U}        
\catcode`\‡=\active\def‡{\'a}        
\catcode`\Ž=\active\defŽ{\'e}        
\catcode`\'=\active\def'{\'{\i}}     
\catcode`\—=\active\def—{\'o}        
\catcode`\œ=\active\defœ{\'u}        
\catcode`\ƒ=\active\defƒ{\'E}        
\catcode`\ˆ=\active\defˆ{\`a}        
\catcode`\=\active\def{\`e}        
\catcode`\"=\active\def"{\`{\i}}     
\catcode`\˜=\active\def˜{\`o}        
\catcode`\=\active\def{\`u}        
\catcode`\Ë=\active\defË{\`A}        
\catcode`\‹=\active\def‹{\~a}        
\catcode`\–=\active\def–{\~n}        
\catcode`\›=\active\def›{\~o}        
\catcode`\Ì=\active\defÌ{\~A}        
\catcode`\"=\active\def"{\~N}        
\catcode`\Í=\active\defÍ{\~O}        
\catcode`\‰=\active\def‰{\^a}        
\catcode`\=\active\def{\^e}        
\catcode`\"=\active\def"{\^{\i}}     
\catcode`\™=\active\def™{\^o}        
\catcode`\ž=\active\defž{\^u}        

\let\optionkeymacros\null

%% file: SSIND.tex
\font\eightrm=cmr8
\font\eighti=cmmi8
\font\eightsy=cmsy8
\font\eightbf=cmbx8
\font\eighttt=cmtt8
\font\eightit=cmti8
\font\eightsl=cmsl8
\font\sixrm=cmr6
\font\sixi=cmmi6
\font\sixsy=cmsy6
\font\sixbf=cmbx6
\skewchar\eighti='177 \skewchar\sixi='177
\skewchar\eightsy='60 \skewchar\sixsy='60
\catcode`\@=11

\def\pc#1#2|{{\bigf@ntpc #1\penalty\@MM\hskip\z@skip\smallf@ntpc #2}}
\def\tenpoint{%
\textfont0=\tenrm \scriptfont0=\sevenrm \scriptscriptfont0=\fiverm
\def\rm{\fam\z@\tenrm}%
\textfont1=\teni \scriptfont1=\seveni \scriptscriptfont1=\fivei
\def\oldstyle{\fam\@ne\teni}%
  \textfont2=\tensy \scriptfont2=\sevensy \scriptscriptfont2=\fivesy
  \textfont\itfam=\tenit
  \def\it{\fam\itfam\tenit}%
  \textfont\slfam=\tensl
  \def\sl{\fam\slfam\tensl}%
  \textfont\bffam=\tenbf \scriptfont\bffam=\sevenbf
  \scriptscriptfont\bffam=\fivebf
  \def\bf{\fam\bffam\tenbf}%
  \textfont\ttfam=\tentt
  \def\tt{\fam\ttfam\tentt}%
  \abovedisplayskip=6pt plus 2pt minus 6pt
  \abovedisplayshortskip=0pt plus 3pt
  \belowdisplayskip=6pt plus 2pt minus 6pt
  \belowdisplayshortskip=7pt plus 3pt minus 4pt
  \smallskipamount=3pt plus 1pt minus 1pt
  \medskipamount=6pt plus 2pt minus 2pt
  \bigskipamount=12pt plus 4pt minus 4pt
  \normalbaselineskip=12pt
  \setbox\strutbox=\hbox{\vrule height8.5pt depth3.5pt width0pt}%
  \let\bigf@ntpc=\tenrm \let\smallf@ntpc=\sevenrm
  \normalbaselines\rm}
\def\eightpoint{%
  \textfont0=\eightrm \scriptfont0=\sixrm \scriptscriptfont0=\fiverm
  \def\rm{\fam\z@\eightrm}%
  \textfont1=\eighti \scriptfont1=\sixi \scriptscriptfont1=\fivei
  \def\oldstyle{\fam\@ne\eighti}%
  \textfont2=\eightsy \scriptfont2=\sixsy \scriptscriptfont2=\fivesy
  \textfont\itfam=\eightit
  \def\it{\fam\itfam\eightit}%
  \textfont\slfam=\eightsl
  \def\sl{\fam\slfam\eightsl}%
  \textfont\bffam=\eightbf \scriptfont\bffam=\sixbf
  \scriptscriptfont\bffam=\fivebf
  \def\bf{\fam\bffam\eightbf}%
  \textfont\ttfam=\eighttt
  \def\tt{\fam\ttfam\eighttt}%
  \abovedisplayskip=9pt plus 2pt minus 6pt
  \abovedisplayshortskip=0pt plus 2pt
  \belowdisplayskip=9pt plus 2pt minus 6pt
  \belowdisplayshortskip=5pt plus 2pt minus 3pt
  \smallskipamount=2pt plus 1pt minus 1pt
  \medskipamount=4pt plus 2pt minus 1pt
  \bigskipamount=9pt plus 3pt minus 3pt
  \normalbaselineskip=9pt
  \setbox\strutbox=\hbox{\vrule height7pt depth2pt width0pt}%
  \let\bigf@ntpc=\eightrm \let\smallf@ntpc=\sixrm
  \normalbaselines\rm}

\newskip\LastSkip
\def\nobreakatskip{\relax\ifhmode\ifdim\lastskip>\z@
  \LastSkip\lastskip\unskip\nobreak\hskip\LastSkip
  \fi\fi}
\catcode`\;=\active
\catcode`\:=\active
\catcode`\!=\active
\catcode`\?=\active
\def;{\nobreakatskip\string;}
\def:{\nobreakatskip\string:}
\def!{\nobreakatskip\string!}
\def?{\nobreakatskip\string?}
%
  \let\titlefont=\seventeenrm

\font\tenss=cmss10 \font\tencaps=cmcsc10
\newfam\ssfam   \textfont\ssfam=\tenss     \def\ss{\fam\ssfam\tenss}
\newfam\capsfam \textfont\capsfam=\tencaps \def\petcap{\fam\capsfam\tencaps}

\font\teneuf=eufm10 \font\seveneuf=eufm7 \font\fiveeuf=eufm5
\font\tenmsa=msam10 \font\sevenmsa=msam7 \font\fivemsa=msam5
\font\tenmsb=msbm10 \font\sevenmsb=msbm7 \font\fivemsb=msbm5

\newfam\msafam \textfont\msafam=\tenmsa \scriptfont\msafam=\sevenmsa
\scriptscriptfont\msafam=\fivemsa
\newfam\msbfam \textfont\msbfam=\tenmsb \scriptfont\msbfam=\sevenmsb
\scriptscriptfont\msbfam=\fivemsb 
\newfam\euffam \textfont\euffam=\teneuf \scriptfont\euffam=\seveneuf
\scriptscriptfont\euffam=\fiveeuf 

\newskip\LastSkip
\def\nobreakatskip{\relax\ifhmode\ifdim\lastskip>0pt
  \LastSkip\lastskip\unskip
  \nobreak\hskip\LastSkip
  \fi\fi}
\catcode`\;=\active \def;{\nobreakatskip\string;}
\catcode`\:=\active \def:{\nobreakatskip\string:}
\catcode`\!=\active \def!{\nobreakatskip\string!}
\catcode`\?=\active \def?{\nobreakatskip\string?}

\frenchspacing
\tenpoint


\magnification 1200

\pretolerance=500   \tolerance=1000   
\brokenpenalty=5000 

\hoffset=1cm
\hsize=125mm \vsize=187mm \parindent=1cm

\parskip=5pt plus 1pt
\newif\ifpagetitre 				\pagetitretrue
\newtoks\hautpagetitre			\hautpagetitre={\hfil}
\newtoks\baspagetitre       \baspagetitre={\hfil\tenrm\folio\hfil}
\newtoks\auteurcourant     \auteurcourant={\hfil}
\newtoks\titrecourant       \titrecourant={\hfil}
\newtoks\hautpagegauche  \newtoks\hautpagedroite
\newtoks\chapitre              \newtoks\paragraphe
\hautpagegauche={{\tenrm\folio}\hfil\the\auteurcourant\hfil}
\hautpagedroite={\hfil\the\titrecourant\hfil{\tenrm\folio}}
\newtoks\baspagegauche     \baspagegauche={\hfil}
\newtoks\baspagedroite     \baspagedroite={\hfil}
\headline={\ifpagetitre\the\hautpagetitre
\else\ifodd\pageno\the\hautpagedroite
\else\the\hautpagegauche\fi\fi}
\footline={\ifpagetitre\the\baspagetitre
\global\pagetitrefalse
\else\ifodd\pageno\the\baspagedroite
\else\the\baspagegauche\fi\fi}

\catcode`\@=11
\def\p@int{{\rm .}}
\def\p@intir{\discretionary{\rm .}{}{\rm .\kern.35em---\kern.7em}}
\def\pointir{\afterassignment\pointir@\let\next=}
\def\pointir@{\ifx\next\par\p@int\else\p@intir\fi\next}
\catcode`\@=12

\def\cite#1{[#1]}
\def\buildo#1\over#2{\mathrel{\mathop{\null#2}\limits^{#1}}}
\def\buildu#1\under#2{\mathrel{\mathop{\null#2}\limits_{#1}}}

\def\_#1{_{\vtop{\halign{$\scriptstyle{##}$&
                       $\scriptstyle{{}##}$\cr #1\crcr}}}}

\def\text#1{\hbox{\rm #1}} 

\def\article#1|#2|#3|#4|#5|#6|#7|
  {\leftskip=7mm\noindent
   \hangindent=2mm\hangafter=1
    \llap{\bf [#1]\hskip1.35em}{\petcap #2}\pointir {\sl #3}, {\rm #4},
    \nobreak {\bf #5} ({\oldstyle #6}), \nobreak #7.\par}
\def\livre#1|#2|#3|#4|#5|
 {\leftskip=7mm\noindent
   \hangindent=2mm\hangafter=1
    \llap{\bf [#1]\hskip1.35em}{\petcap #2}\pointir {\sl #3}, {\rm #4},
    {\oldstyle #5}.\par}
\def\alivre#1|#2|#3|#4|#5|#6|#7|
  {\leftskip=7mm\noindent
    \hangindent=2mm\hangafter=1
    \llap{\bf [#1]\hskip1.35em}{\petcap #2}\pointir {\sl #3},
dans {\it #4}, \'edit\'e par {\bf #5}, {\rm #6}, {\oldstyle #7}.\par}
\def\divers#1|#2|#3|#4|
  {\leftskip=7mm\noindent
   \hangindent=2mm\hangafter=1
    \llap{[#1]\hskip1.35em}{\petcap #2}\pointir #3,
    {\oldstyle #4}.\par}
\newcount\secno \secno=0
\newcount\ssecno \ssecno=0
\newcount\sssecno \sssecno=0
\newcount\chapno \chapno=0
\newcount\notenumber \notenumber=1
\newcount\exino \exino=0
\newcount\expno \expno=0
%
%
%
\newdimen\indentsec\indentsec=20pt
\newdimen\indentssec\indentssec=20pt
\newdimen\indentsssec\indentsssec=20pt
\newdimen\indentrem\indentssec=0pt

%
\newdimen\indentTh\indentTh=0pt

%
\newdimen\indentth\indentssec=0pt
\def\titregen#1#2|{\par\vskip .5cm\penalty -100
             {\leftskip=0pt plus \hsize
			   \rightskip=\leftskip
			   \parfillskip=0pt
			   \baselineskip=17pt
			   \noindent
             \titlefont{ #1}{#2}\par}
			  \vskip 5pt\penalty 500}
\def\titre#1|{\titregen{}{#1}|}
%
\def\ntitre#1|{\global\advance\chapno by 1\global\ssecno=0\global\secno=0
\titregen{\the\chapno\ }{#1}|}
\def\auteur#1|{\vskip 5pt\penalty 100
               \vbox{\centerline{\sl par #1}
                    \vskip 5pt}\penalty 100}
\def\sectiongen#1#2#3{\parindent=\indentsec\par\vskip .3cm
\vskip 0mm plus -20mm minus 1,5mm\penalty-50
{\bf #1}{\bf #2}{#3}\nobreak\parindent=20pt}%
%
\def\secc#1|{\sectiongen{}{#1}{\pointir}}
%
\def\nsecc#1|%
{\global\advance\secno by 1\global\ssecno=0\global\sssecno=0
\sectiongen{\the\secno\ }{#1}{\pointir}}
%
%
\def\secp#1|{\sectiongen{}{#1}{}\par}
%
\def\nsecp#1|%
{\global\advance\secno by 1\global\ssecno=0\global\sssecno=0
\sectiongen{\the\secno\ }{#1}{}\par}
\def\ssectiongen#1#2#3#4{\parindent=\indentssec\par\vskip .2cm
\vskip 0mm plus -20mm minus 1,5mm\penalty-50
{\bf #1}{\sl #2}{\sl #3}{#4}\nobreak\medskip\parindent=20pt}%
%
\def\ssecc#1|#2{\ssectiongen{}{#1}{\pointir}{#2}}
%
\def\nssecc#1|#2{\global\advance\ssecno by 1\global\sssecno=0
\ssectiongen{\the\secno.\the\ssecno\ }{#1}{\pointir}{#2}}
%
\def\ssecp#1|{\ssectiongen{}{#1}{}{}\par}
%
\def\nssecp#1|{\global\advance\ssecno by 1\global\sssecno=0
\ssectiongen{\the\secno.\the\ssecno\ }{#1}{}{}\par}
\def\sssectiongen#1#2#3#4{\parindent=\indentsssec\par\vskip .2cm
\vskip 0mm plus -20mm minus 1,5mm\penalty-50
{\bf #1}{\sl #2}{\sl #3}{#4}\nobreak\medskip\parindent=20pt}%
%
\def\sssecc#1|#2{\sssectiongen{}{#1}{\pointir}{#2}}
%
\def\nsssecc#1|#2{\global\advance\sssecno by 1
\ssectiongen{\the\secno.\the\ssecno.\the\sssecno\ }{#1}{\pointir}{#2}}
%
\def\sssecp#1|{\sssectiongen{}{#1}{}{}\par}
%
\def\nsssecp#1|{\global\advance\sssecno by 1
\sssectiongen{\the\secno.\the\ssecno.\the\sssecno\ }{#1}{}{}\par}
\long\def\Th#1#2#3#4{\parindent=\indentTh\par\vskip5pt
{#1}{\petcap #2}{\sl #3}\parindent=20pt{\sl #4\par}\vskip 5pt\parindent=20pt}
\long\def\th#1#2#3#4{\parindent=\indentth\par\vskip5pt
{#1}{\pppetcap #2}{\eightpoint\sl #3}\parindent=20pt{\eightpoint\sl #4\par}\vskip 5pt\parindent=20pt}
\long\def\remarque#1#2#3#4{\parindent=\indentrem\par\vskip5pt
{#1}{\eightpoint \sl #2}{\eightpoint\sl #3}\parindent=20pt{\eightpoint#4\par}
\vskip 5pt\parindent=20pt}
\long\def\remarques#1#2#3#4{\parindent=\indentrem\par\vskip5pt
{#1}{\eightpoint \sl #2}{\eightpoint\sl #3}{\eightpoint#4\par}
\vskip 5pt\parindent=20pt}

\long\def\remarquesa#1#2#3#4#5{\parindent=\indentrem\par\vskip5pt
{#1}{\eightpoint \sl #2}{\eightpoint\sl #3}{\eightpoint#4}
\parindent=20pt{\eightpoint#5\par}
\vskip 5pt\parindent=20pt}

\long\def\Remarque#1#2#3#4{\parindent=\indentrem\par\vskip5pt
{#1}{ \sl #2}{\sl #3}\parindent=20pt{#4\par}
\vskip 5pt\parindent=20pt}
\long\def\remarquesn#1#2#3#4{\parindent=\indentrem\par\vskip5pt
{\eightpoint \sl #1}{\eightpoint #2}{\eightpoint\sl #3}{\eightpoint#4\par}
\vskip 5pt\parindent=20pt}

\long\def\Remarquen#1#2#3#4{\parindent=\indentrem\par\vskip5pt
{ \sl #1}{#2}{\sl #3}\parindent=20pt{#4\par}
\vskip 5pt\parindent=20pt}

%
\long\def\Thc#1|#2\finc{\Th{}{#1}{\pointir}{#2}}
%
\long\def\Thnc#1|#2|#3\finnc{\Th{#1}{#2}{\pointir}{#3}}
\long\def\Exic#1|#2\finc{{\global\advance\exino by 1}\Remarquen%
{Exercice }{\the\exino}{ #1\pointir}{#2}}
\long\def\Expc#1|#2\finc{{\global\advance\expno by 1}\Remarquen%
{Exemple }{\the\expno}{ #1\pointir}{#2}}
\long\def\exic#1|#2\finc{{\global\advance\exino by 1}\remarquesn%
{\bf Exercice }{\the\exino}{ #1\pointir}{#2}}
\long\def\expc#1|#2\finc{{\global\advance\expno by 1}\remarquesn%
{Exemple }{\the\expno}{ #1\pointir}{#2}}

%
\long\def\Ec#1\finc{\Th{}{}{}{#1}}
%

\long\def\thc#1|#2\finc{\th{}{#1}{\pointir}{#2}}%
\long\def\thc#1\finc{\th{}{}{}{#1}}

\long\def\Thp#1|#2\finp{\Th{}{#1}{\par}{#2}}
%
\long\def\thp#1|#2\finp{\th{}{#1}{\par}{#2}}
%
\long\def\rmc#1|#2\finc{\remarque{}{#1}{\pointir}{#2}}
%
\long\def\Rmc#1|#2\finc{\Remarque{}{#1}{\pointir}{#2}}

\long\def\rmp#1|#2\finp{\remarque{}{#1}{\par}{#2}}
%
\long\def\Rmp#1|#2\finp{\Remarque{}{#1}{\par}{#2}}

\long\def\parc#1\finc{\remarque{}{}{}{#1}}
%
\long\def\parcs#1\fincs{\remarques{}{}{}{#1}}

\long\def\parcsa#1\fins#2\fincsa{\remarquesa{}{}{}{#1}{#2}}

\def\Rm#1|{\parindent=0pt\par\vskip5pt{\sl #1}\pointir\parindent=20pt}
%

%


%

\def\preuved#1|{\parindent=0pt\par{\sl Preuve d#1}\pointir\parindent=20pt}

\def\qed{\quad\hbox{\hskip 1pt\vrule width 4pt height 6pt
         depth 1.5pt\hskip 1pt}}
\def\findem{\penalty 500 \hbox{\qed}\par\vskip 3pt}

%

%

%

%

%
\def\og{\leavevmode\raise.3ex\hbox{$\scriptscriptstyle\langle\!\langle\,$}}

\def\fgf{\/\leavevmode\raise.3ex\hbox{$\scriptscriptstyle\,\rangle\!\rangle$}}
%
\def\fg{\fgf\ }
\def\note#1{\footnote{$^{\the\notenumber}$}{#1}%
\global\advance\notenumber by 1}%
\def\ieme{\raise 1ex\hbox{\pc{}i\`eme|}\ }
\def\iemes{\raise 1ex\hbox{\pc{}i\`emes|}\ }

%% file: FONTESMATH.tex
\input amsfonts.def

%% file: MATH
\def\build#1_#2^#3{\mathrel{\mathop{\kern 0pt#1}\limits_{#2}^{#3}}}
%

\def\hdfl#1#2{\smash{\mathop{\hbox to 12mm{\rightarrowfill}}
\limits^{\scriptstyle#1}_{\scriptstyle#2}}}

\def\hdhfl#1#2{\smash{\mathop{\hbox to 12mm{\hookrightarrowfill}}
\limits^{\scriptstyle#1}_{\scriptstyle#2}}}

\def\hgfl#1#2{\smash{\mathop{\hbox to 12mm{\leftarrowfill}}
\limits^{\scriptstyle#1}_{\scriptstyle#2}}}

\def\hghfl#1#2{\smash{\mathop{\hbox to 12mm{\hookleftarrowfill}}
\limits^{\scriptstyle#1}_{\scriptstyle#2}}}













\def\cotg{\mathop{\rm cotg}\nolimits}

\def\th{\mathop{\rm th}\nolimits}


%% file: pdffix.tex
\input ifpdf.sty
\ifpdf
\pdfpagewidth=8.5truein
\pdfpageheight=11truein
\pdfhorigin=1truein
\pdfvorigin=1truein
\fi

%% file: C0.tex
%

\auteurcourant={\hfill   $\sum{{\sin(kx)}\over{k}}$\hfill {}}
\titrecourant={\hfill\eightrm Capo des s\'eries de Fourier\hfill}
 \null\vskip-15mm
\parc
\null\vskip-5mm
\titre  Le capo%
\note{\eightpoint
chef en italien,  dans le
Gr\'esivaudan  la vache dominante conduisant le
troupeau \`a l'alpage.\hfill
} du troupeau des s\'eries de Fourier.|
\finc
\vskip-3mm
\auteur{
$\sum{{\sin(kx)}\over{k}}$}|
\vskip 2mm
{\eightpoint
\centerline{ Valence-Lyon-Grenoble 
printemps automne hiver {\oldstyle 2011}, {\oldstyle 2012}}

 %
%
\centerline{Collection 
de \og cours {\`a} distance\fg
d'accoutumance {\`a}
de robustes\note{\eightpoint
mais que sous le r{\`e}gne de El'hemd{\'e}e, certains 
traitent de poussi{\`e}reuses.
 }
r{\'e}f{\'e}rences.
}
}
\vskip2mm
\parc
\centerline{\bf R\'esum\'e}
$\displaystyle
x-\pi+{{\sin(nx)}\over{n}}+2\sum_{k=1}^{n-1}{{\sin(kx)}\over{k}}=
\int_{\pi}^x[1+\cos(nt)+2\sum_{k=1}^{n-1}\cos(kt)]dt=
\int_{\pi}^x\sin(nt)\cotg({{t}\over{2}}) dt$

Par int\'egration par parties,
$\displaystyle
S_n^{\ast}(x)\!=\!{{\sin(nx)}\over{2n}}+\!\sum_{k=1}^{n-1}\!\!{{\sin(kx)}\over{k}}$,
 les sommes partielles modifi\'ees%
\note{\eightpoint
qui donnent lieu \`a des noyaux de Dirichlet
d'\'ecriture et majoration plus simple que sommes et noyaux de Dirichlet usuels.
Avec, dans leurs preuves, 
 calligraphie et majorations plus lourdes,
les \'enonc\'es A, B et C ci-dessous ont aussi lieu pour  sommes
et noyaux de Dirichlet usuels.},
sont dans 
$]\!-{{\pi}\over{2}}, {{\pi}\over{2}}[$ et,
pour tout $\epsilon\!\in\!]0, \pi[$,
elles convergent, quand $n$ tend vers l'infini, unifor\-m\'e\-ment 
 sur $]\epsilon, 2\pi-\epsilon[$ vers
${{\pi-x}\over{2}}$.
Donc la {\it corne\/}
$2\pi$-p\'eriodique
$\psi : {\bf R\/}\rightarrow{\bf R\/}, \psi(0)\!=\!0, 0\!<\!x<\!2\pi, \psi(x)\!=\!{{\pi-x}\over{2}}$
est limite en moyenne quadratique de sa s\'erie de Fourier et~:

$A$  Si $0<\alpha<\beta<2\pi$ alors
$\chi_{\alpha, \beta} : [0, 2\pi[\rightarrow{\bf R},
\chi_{\alpha, \beta}(x)\!=\!{{1}\over{\pi}}\bigl[\psi(x-\alpha)-{{\alpha}\over{2}}
+{{\beta}\over{2}}-\psi(x-\beta)\bigr]$
est la fonction caract\'eristique\note{\eightpoint
modifi\'ee aux extr\'emit\'ees par  $\chi(\alpha)={{1}\over{2}}=\chi(\beta)$.}
$\chi_{[\alpha, \beta]}$ de l'intervalle $[\alpha, \beta]$.
Ainsi toute fonction en escalier est limite  en moyenne quadratique 
de sa s\'erie de Fourier  et,
par densit\'e des fonctions en escalier dans
les fonctions de carr\'e int\'egrable, le m\^eme r\'esultat 
a lieu pour toute 
fonction de carr\'e int\'egrable.

$B$ {\sl La primitive 
$E^{\ast}_n(x)\!=\!\int_{\pi}^x\!\!D^{\ast}_n$
de
 $\displaystyle D^{\ast}_n(y)\!=\!1\!+\!\cos(ny)\!+\!%
2\!\!\sum_{k=1}^{n-1}\!\!\cos(ky)$,
le {\it  $n$-noyau de Dirichlet modifi\'e\/},
est
 \`a valeurs dans 
$[-\pi, \pi]$
 et, sur tout 
$[\epsilon, 2\pi-\epsilon], 0\!<\!\epsilon\!<\!\pi$,
tend uniform\'ement vers $0$.
}

$C$ Les sommes de Fourier partielles modifi\'ees  d'une fonction 
$2\pi$ p\'eriodique int\'egrable sont :
$$\displaystyle
 S^{\ast}_n(f)(x)\!=\!{{1}\over{2}}%
(c_{-n}(f)e^{-inx}\!\!+\!c_{n}(f)e^{inx})
\!+\!\!\!%
\sum_{k=1-n}^{n-1}\!\!\!c_k(f)e^{ikx}\!=\!%
{{1}\over{2\pi}}\int_x^{x+2\pi}\!\!\!\!\!f(t)D^{\ast}_n(t-x)dt$$
D'o\`u, par int\'egration par parties sur chaque intervalle de d\'erivabilit\'e,
le th\'eor\`eme  de Dirichlet :

{\sl  Si $f$ est $2\pi$-p\'eriodique et, pour
$x\!=\!\alpha_0\!<\!\cdots\!<\!\alpha_{m+1}\!=\!x+2\pi$,
  d\'erivable sur chaque intervalle ouvert
$]\alpha_i, \alpha_{i+1}[, i\!=\!0,\ldots,m$
avec $f'$ int\'egrable alors pour tout 
$x\!\in\!{\bf R}$
on a~:
$S^{\ast}_n(f)(x)\!=\!$
$${{1}\over{2\pi}}\Bigl[f_{-}(x)E^{\ast}_n(2\pi)-f_{+}(x)E^{\ast}_n(0)%
-\bigl[\sum_{j=1}^{m}[f_{+}(\alpha_j)-%
f_{-}(\alpha_j)]E^{\ast}_n(\alpha_j)+%
\int_x^{x+2\pi}\!\!\!\!\!\!f'(t)E^{\ast}_m(t-x)dt\bigr]\Bigr]$$
tend,  quand
$n\!\rightarrow\!\infty$, vers la moyenne
  $\displaystyle {{f_{-}(x)+\!f_{+}(x)}\over{2}}$
des limites \`a gauche et \`a droite de $f$ en $x$.
}
\finc
\vskip1mm
\parc
\centerline{\bf Abstract}
A glimpse at 
$\sum{{\sin(kx)}\over{k}}$ gives  the above few
  lines exposition of Fourier series's 
quadratic mean convergence
for square integrable functions and Dirichlet's conver\-gence
theorem of the Fourier serie of a piecewise differentiable function
with integrable derivative.
\finc


%% file: C1.tex
\null\vskip1mm
{\eightpoint
\centerline{\eightpoint\petcap
 Commentaires bibliographiques\/}
\vskip7mm

Les premiers textes ({\bf [D]\/}\note{\eightpoint
reproduit dans {\bf [K L-R]\/} avec une histoire
d\'etaill\'ee des s\'eries de Fourier cf. aussi
{\bf [Le]\/}\S{\bf 15\/}-{\bf 19\/}.%
}%
, chapitre IV du  tome II de {\bf [J2]\/}, 
 chapitre III de {\bf [Le]\/}, Kapitel 31 de {\bf [La]\/}\dots)
ignorent\note{\eightpoint
cependant d\`es 1903 Hurewitz {\bf [H]\/} met l'accent sur la th\'eorie
quadratique en la d\'eduisant du r\'esultat, que F\'ejer venait d'\'etablir,
de convergence uniforme vers une fonction continue des moyennes des
$n$ premi\`eres sommes partielles de sa s\'erie de Fourier, cf. aussi
{\bf [Z]\/}III \S3 p. 89.
}
la th\'eorie quadratique,  commen\c{c}ant par
 le th\'eor\`eme de Dirichlet pour les fonctions monotones par morceaux
et, apr\`es Jordan {\bf [J1]\/}, \`a variation born\'ee.

Ce r\'esultat plus g\'en\'eral
que C 
se montre cependant comme ici, une fois
not\'e que l'int\'egration (de Stieljes) par parties a 
lieu 
pour les fonctions \`a variation born\'ee
 (cf. \S{\bf 54\/}  de {\bf [R Sz-N]\/}).

La compl\'etude du syst\`eme trigonom\'etrique%
\note{\eightpoint  d'o\`u
I {\sl La s\'erie de Fourier de $f$ de carr\'e int\'egrable
converge en moyenne quadratique vers $f$.\/} et l'\'enonc\'e ponctuel
II {\sl Si les sommes partielles d'une s\'erie trigonom\'etrique $\Sigma$ sont
uniform\'ement major\'ees et convergent simplement vers une fonction $f$
alors $\Sigma$ est la s\'erie de Fourier de $f$.}
}  aujourd'hui
({\bf 2.6\/} de {\bf [H R]\/}, {\bf 4.24\/} de {\bf [Ru]\/}\dots)
suit souvent de la densit\'e 
des poly\-n\^omes trigonom\'etriques dans le Banach des
fonctions continues.
 
Lebesgue ({\bf [Le]\/}\S{\bf 24\/}, {\bf [Z]\/}I\S6.)
obtenait la compl\'etude  directement\note{\eightpoint
Par contrapos\'ee, {\sl $f$  continue $2\pi$-p\'eriodique  non nulle
a un  coefficient de Fourier non nul\/}~:
Soit $\alpha\!<\!\beta\!<\!\alpha+\pi$ et $\epsilon\!>\!0$
tels que que pour $x\!\in\![\alpha, \beta], |f(x)|\!\geq\!\epsilon$,
quitte \`a prendre $-f, f(x)\!\geq\!\epsilon$.
Comme 
$\cos(x-{{\alpha+\beta}\over{2}})-\cos({{\beta-\alpha}\over{2}})\in%
]-2, 1]$,
est nulle en $\alpha, \beta$ 
et, entre  ${{\alpha+\beta}\over{2}}$
et ${{\alpha+\beta}\over{2}}\pm\pi$,
 monotone , le polyn\^ome trigonom\'etrique
$\displaystyle\psi_n\!=\!\Bigl[1+{{e^{i({{\alpha+\beta}\over{2}}-x)}%
+e^{i(x-{{\alpha+\beta}\over{2}})}}\over{2}}-%
\cos({{\beta-\alpha}\over{2}})\Bigr]^n\!=\!%
\sum_{k=-n}^{n}d_ke^{-ikx}$
tend vers $+\infty$ sur $]\alpha, \beta[$ et $0$ sur
$]\beta, \alpha+2\pi[$ donc pour $n$ assez grand~:\hfill\break
$\displaystyle\int_0^{2\pi}f\psi_n\!=\!2\pi\!\!\sum_{k=-n}^nd_kc_k(f)\!>\!0$
  et un des coefficients de Fourier, $c_k(f)$ 
de $f$ est non nul.%
}, 
sans utiliser la plus fine densit\'e uniforme, citant
une preuve simultan\'ee\note{\eightpoint
trop longue pour \^etre expos\'ee ici, mais avec des techniques Hilbertiennes 
\'etonnamment alg\'ebriques et un argument de s\'erie enti\`ere,
elle semble oubli\'ee (tant par {\bf [K L-R]\/} que {\bf [Z]\/}).%
} {\bf [K]\/}I\S{\bf 2\/} d'A. Kneser.

Tous deux  ({\bf [Le]\/}\S{\bf 1\/}, \S{\bf 28\/},
{\bf [K]\/}I\S{\bf 4\/}-{\bf 5\/})
donnent aussi
 une r\'eduction \`a la compl\'e\-tude
du th\'eor\`eme de Dirichlet pour les 
fonctions $C^2$ par morceaux\note{\eightpoint
utilisant le second th\'eor\`eme de la moyenne,
Lebesgue demande seulement \`a la d\'eriv\'ee d'\^etre \`a
variation born\'ee, en fait cette preuve passe pour $f$ r\'egl\'ee
absolument continue sur ses intervalles de continuit\'e,
  \'enonc\'e que Hardy et Rogosinski ({\bf [H R]\/} \S {\bf 3.8\/}-{\bf 3.9\/})
font suivre de la compl\'etude (note $^3$ ci-dessus) 
et d'une int\'egration par partie
(\S {\bf 3.4\/}) pour chaque coefficient de Fourier.%
}
\`a nombre fini de discontinuit\'es, le r\^ole de
$
\sum\!\!{{\sin(kx)}\over{k}}$ \'etant plus explicite 
dans {\bf [Le]\/} que dans 
{\bf [K]\/}~:

{\sl Soit $f$ $2\pi$-p\'eriodique  \`a variation born\'ee 
(resp. d\'erivable par morceaux) alors~:

1) $f$ est somme d'une fonction continue 
\`a variation born\'ee (resp. d\'erivable par morceaux) et
du mod\`ele explicite 
$\displaystyle\sum_{0\leq \alpha<2\pi}{{f_{+}(\alpha)-f_{-}(\alpha)}\over{\pi}}%
\sum_{k=1}^{\infty}{{\sin(k(x-\alpha))}\over{k}}$.

2) Si  $f$ continue et $f'$ 
\`a variation born\'ee%
\note{\eightpoint Une remarque de de la Vall\'ee Poussin
({\bf [B]\/} \S{\bf 5},
{\bf [dVP]\/} p. 32),
\'etend 2) quand
$f'$
est seulement de carr\'e int\'egrable~: Les in\'egalit\'es de Bessel 
$\displaystyle  2\pi\sum|c_n(f')|^2\!\leq\!\int_0^{2\pi}\!\!|f'|^2$,
de Cauchy-Schwarz et
 $c_n(f)\!=\!{{1}\over{in}}c_n(f')$ donnent
 $\Sigma c_n(f)$ absolument convergente~: II de la note
$^3$ p. 2. s'applique} 
alors $f$ est limite uniforme de $S^{\ast}_n(f)$.}

Hardy\&Rogosinski enfin, apr\`es l'avoir d\'eduite de la compl\'etude et d'une
majo\-ration de ses sommes partielles\note{\eightpoint
$|S_n^{\ast}(x)|\!\leq\!%
{\sum_{k={{1}\over{2}}}^{{\pi}\over{2|x|}}}^{\!\!\!\ast}\!{{|\sin(kx)|}\over{k}}+%
\chi_{[0, n]}({{\pi}\over{2|x|}})%
|{\sum_{k={{\pi}\over{2|x|}}}^n}^{\!\!\!\!\!\!\ast}\ {{\sin(kx)}\over{k}}|\!\leq\!{{\pi}\over{2}}+%
{{2|x|}\over{\pi}}\cotg({{x}\over{2}})\!\leq\!{{\pi}\over{2}}+{{4}\over{\pi}}$.%
}
 obtiennent aussi ({\bf [H R]\/} fin du \S{\bf 3.7\/})
la somme de  $\displaystyle\sum{{\sin(kx)}\over{k}}$
 \`a partir, 
suivant
Euler, du d\'eveloppement\note{\eightpoint
et du th\'eor\`eme d'Abel en $z_0\!=\!e^{inx}$,
justifiant l'intuition d'Euler,
cf. aussi {\bf [Z]\/}I\S2, {\bf [Le]\/}\S{\bf 21\/}.%
}, pour $|z|\!<\!1$, de
$\log(1-z)$.

Int\'egrer par partie la d\'eriv\'ee qui se somme explicitement
pour d\'eterminer une somme (comme ici le {\it capo\/}
$\displaystyle \sum{{\sin(kx)}\over{k}}$) remonte \`a Fourier
({\bf [F]\/} Art. 178-180).

}

%% file: Ref.tex
\null
{\eightpoint
\vskip3mm
\centerline{R\'EF\'ERENCES}
\vskip3mm

{\hangindent=1cm\hangafter=1\noindent{\bf   [B]}\
{\petcap  M. B\^ocher}\pointir
 {\sl Introduction to the theory of Fourier's series},\hfill\break
Annals of Mathematics, 7, 81-152,
{\oldstyle 1906}.
\par}

{\hangindent=1cm\hangafter=1\noindent{\bf   [D]}\
{\petcap  P.J. L. Dirichlet}\pointir
 {\sl Sur la convergence des s\'eries trigonom\'etriques qui servent
\`a repr\'esenter une fonction arbitraire entre deux limites donn\'ees},\hfill\break
Journal f\H ur die reine und angewandte Mathematik (Journal de Crelle), 4, 157-169,
{\oldstyle 1829}.
\par}

{\hangindent=1cm\hangafter=1\noindent{\bf   [F]}\
{\petcap  J. Fourier}\pointir
 {\sl Th\'eorie analytique de la chaleur},
Firmin Didot, p\`ere et fils,
{\oldstyle 1822}.
\par}

{\hangindent=1cm\hangafter=1\noindent{\bf   [H R]}\
{\petcap  G.H. Hardy, W.W. Rogosinski}\pointir
 {\sl Fourier series},
Cambridge tracts in mathematics and mathematical physics 38,
Cambridge University Press
{\oldstyle 1956}.
\par}

{\hangindent=1cm\hangafter=1\noindent{\bf   [H]}\
{\petcap  A. Hurwitz}\pointir
 {\sl \"Uber die Fourierschen Konstanten integriebarer Funktionen},\hfill\break
Mathematische Annalen, 57, 425-446
{\oldstyle 1903}.
\par}

{\hangindent=1cm\hangafter=1\noindent{\bf   [J1]}\
{\petcap  C. Jordan}\pointir
 {\sl Sur la s\'erie de Fourier},\hfill\break
Comptes rendus de l'Acad\'emie des sciences, Paris, 92, 228-230,
{\oldstyle 1881}.
\par}

{\hangindent=1cm\hangafter=1\noindent{\bf   [J2]}\
{\petcap  C. Jordan}\pointir
 {\sl Cours d'analyse de l'\'ecole polytechnique t. I, II, III},\hfill\break
Paris, Gauthier-Villars,
{\oldstyle 1882}.
\par}

{\hangindent=1cm\hangafter=1\noindent{\bf   [K L-R]}\
{\petcap  J.-P. Kahane, P.-G. Lemari\'e-Rieusset}\pointir
 {\sl S\'eries de Fourier et ondelettes},\hfill\break
Paris, Cassini,
{\oldstyle 1998}.
\par} 

{\hangindent=1cm\hangafter=1\noindent{\bf   [K]}\
{\petcap  A. Kneser}\pointir
 {\sl Untersuchungen \H uber die Darstellung willk\H urlicher Funktionen
in der mathematischen Physik},
Mathematische Annalen, 58, 81-147,
{\oldstyle 1904}.
\par}

{\hangindent=1cm\hangafter=1\noindent{\bf   [La]}\
{\petcap  E. Landau}\pointir
 {\sl Einf\H uhrung in die Differentialrechnung und  Integralrechnung},\hfill\break
Groningen, P. Nordhooff,
{\oldstyle 1934}.
\par} 

{\hangindent=1cm\hangafter=1\noindent{\bf   [Le]}\
{\petcap  H. Lebesgue}\pointir
 {\sl Le\c cons sur les s\'eries trigonom\'etriques},\hfill\break
Paris, Gauthier-Villars,
{\oldstyle 1906}.
\par} 

 {\hangindent=1cm\hangafter=1\noindent{\bf   [R, Sz-N]}\
{\petcap  F. Riesz, B. Sz-Nagy}\pointir
 {\sl Le\c cons d'analyse fonctionnelle},\hfill\break
Paris, Gauthier-Villars; Budapest, Akad\'emiai Kiado,
{\oldstyle 1968}.
\par}

{\hangindent=1cm\hangafter=1\noindent{\bf   [Ru]}\
{\petcap  W. Rudin}\pointir
 {\sl Real and complex analysis.} $2^{nd}$ edition,
McGraw-Hill,
{\oldstyle 1974}.
\par} 

{\hangindent=1cm\hangafter=1\noindent{\bf   [VP]}\
{\petcap  C. {\rm de la\/} Vall\'ee Poussin}\pointir
 {\sl Sur quelques applications de l'int\'egrale de Poisson},
Annales de la soci\'et\'e scientifique de Bruxelles 17 B, 18-34,
{\oldstyle 1893}.
\par} 

{\hangindent=1cm\hangafter=1\noindent{\bf   [Z]}\
{\petcap  A. Zygmund}\pointir
 {\sl Trigonometric series} $2^{nd}$ edition,\hfill\break
Cambridge University Press,
{\oldstyle 1959}\quad [$1^{\hbox{\sevenrm i\`ere}}$ \'edition,
Varsovie, {\oldstyle 1935}].
\par} 
}

%% file: C2.tex
\null\vskip-5mm
\secno=-1
\centerline{\petcap Appendices}
\nsecp Sommes partielles et sommes partielles modifi\'ees.|
Les  {\it sommes modifi\'ees\/} 
d'une suite\note{\eightpoint
dans un groupe ab\'elien $A$ o\`u l'on peut diviser par $2$ et
indic\'ees par une partie discr\`ete ${\bf X\/}\!\subset\!{\bf R\/}$
des r\'eels,
ici $A={\bf C\/}, {\bf X\/}\!=\!{\bf N\/}\!\setminus\!\{0\}, {\bf Z\/}$
les entiers positifs ou relatifs
et, si
 $k\!\in\!{\bf R\/}\setminus\!{\bf X\/}, u_k\!=\!0$.
} 
$u_x, x\!\in\!{\bf X\/}$ entre  $a, b\!\in\!{\bf R}$
sont d\'efinies pour   v\'erifier la relation de Chasles
$\displaystyle{\sum_a^c}^{_{\ast}}u_n\!=\!{\sum_a^b}^{_{\ast}}u_n+%
{\sum_b^c}^{_{\ast}}u_n$~: si $a\leq b$
on a~:
$$\displaystyle {\sum_a^b}^{_{\ast}}u_n\!=\!\sum_{a\leq k\leq b}^{}u_k-%
{{1}\over{2}}[u_a+u_b]\!=\!\sum_{a<k\leq b}{{1}\over{2}}(u_{k-1}+u_k)$$
\vskip-3mm
\nsecp L'expression 
$D^{\ast}_n(t)=\sin(nt)\cotg({{t}\over{2}})$
du noyau de Dirichlet modifi\'e.|
La relation
$\sin({{t}\over{2}})D^{\ast}_n(t)=\sin(nt)\cos({{t}\over{2}})$ 
suit, si
 $c_x\!=\!\cos(kt), s_x\!=\!\sin(xt)$%
, de~:\hfill\break
$\displaystyle 
s_{{1}\over{2}}\sum_{k=1}^n(c_{k-1}+c_k)
\!=\!%
s_{{1}\over{2}}\sum_{k=1}^n2c_{k-{{1}\over{2}}}c_{{1}\over{2}}\!=\!%
\bigr[\sum_{k=1}^n2s_{{1}\over{2}}c_{k-{{1}\over{2}}}\bigl]c_{{1}\over{2}}\!=\!%
\bigr[\sum_{k=1}^n(-s_{k-1}+s_{k})\bigl]c_{{1}\over{2}}$
.\hfill\findem



\nsecp La premi\`ere int\'egration par parties.|
$E^{\ast}_n(x)=\int_{\pi}^x\sin(nt)\cotg({{t}\over{2}})dt=
-{{\cos(nx)}\over{n}}\cotg({{x}\over{2}})+
\int_{\pi}^x{{\cos(nt)}\over{n}}{{d}\over{dt}}\cotg({{t}\over{2}}) dt$
donc, si $0<x<2\pi$ on a
$|E^{\ast}_n(x)|\leq{{2}\over{n}}|\cotg({{x}\over{2}})|\leq
{{4}\over{n\min(x, 2\pi-x)}}$
et, si $0\!<\!\epsilon\!<\! \pi$,
 $E^{\ast}_n(x)$ converge vers z\'ero uniform\'ement sur l'intervalle
$[\epsilon, 2\pi-\epsilon]$.

\nsecp La majoration uniforme
$\displaystyle \big{|}E^{\ast}_n(x)\big{|}\!=\!%
\big{|}\int_{\pi}^xD^{\ast}_n(t)dt\big{|}\!\leq\!\pi$.|
 $\displaystyle D^{\ast}_n(x)\!=\!1+\cos(nx)+2\!\sum_{k=1}^{n-1}\!\cos(kx)$
donc $E^{\ast}_n(2\pi-x)\!=\!-E^{\ast}_n(x)$ et
$\displaystyle E^{\ast}_n(0)\!=\!%
-\int_0^{\pi}D^{\ast}_n\!=\!-\pi$,
 il suffit   d'\'etablir que si 
$x\!\in\![0, \pi]$
alors $-\pi\!=\!E^{\ast}_n(0)\!\leq\!E^{\ast}_n(x)\!\leq\!\pi$.
Comme $]0, \pi[\ni t\mapsto\cotg({{t}\over{2}})$ est positive d\'ecroissante,
$E^{\ast}_n$ est, sur $[0, \pi]$, major\'ee par son premier maximum relatif,
qui puisque
${E^{\ast}_n}'(x)\!=\!\sin(nx)\cotg({{x}\over{2}})$ est
$E^{\ast}_n({{\pi}\over{n}})$
 qui, d'apr\`es {\bf 2\/} ci-dessus,
admet la majoration~:
$\displaystyle |E^{\ast}_n({{\pi}\over{n}})|\leq
{{4}\over{n{{\pi}\over{n}}}}\!=\!{{4}\over{\pi}}\!\leq\!\pi$.\hfill\findem

\nsecp Fonctions cr\'eneaux en deux coups de la corne $\psi$.|
Soit $C_{\alpha, \beta}\!=\!2\pi\chi_{\alpha, \beta}$.
Si $x\in]0, 2\pi[\setminus[\alpha, \beta]$ alors
$x-\alpha, x-\beta\!\in]-k2\pi, (1-k)2\pi[$ pour un m\^eme
$k\!=\!1, 0$
 et
$\displaystyle C_{\alpha, \beta}(x)\!=\!%
\pi-(x-\alpha+2k\pi)-\alpha+\beta%
-[\pi-(x-\beta+2k\pi)]\!=\!0$.

Si $x\in]\alpha, \beta[$ alors
$\displaystyle C_{\alpha, \beta}(x)\!=\!%
\pi-(x-\alpha)-\alpha%
+\beta-[\pi-(x-\beta+2\pi)]\!=\!2\pi$.

$\displaystyle C_{\alpha, \beta}(\alpha)\!=\!%
0-\alpha%
+\beta-[\pi-(\alpha-\beta+2\pi)]\!=\!\pi$
et
$\displaystyle C_{\alpha, \beta}(\beta)\!=\!%
\pi-(\beta-\alpha)-\alpha%
+\beta-0\!=\!\pi$.

\nsecp L'int\'egrale de Dirichlet $\displaystyle S^{\ast}_n(f)(x)\!=\!%
{{1}\over{2\pi}}\int_x^{x+2\pi}\!f(t)D^{\ast}_n(t-x)dt$.|
Le $k^{\hbox{\sevenrm i\`eme}}$ coefficient de Fourier de $f$ est
$c_k(f)\!=\!{{1}\over{2\pi}}\int_0^{2\pi}f(t)e^{-ikt}dt$,
donc la $n^{\hbox{\sevenrm i\`eme}}$ som\-me
 modifi\'ee
 $\displaystyle S^{\ast}_n(f)(x)\!=\!%
{{1}\over{2}}\bigl[c_{-n}(f)e^{-inx}+ c_{n}(f)e^{inx}\bigr]+%
\!\sum_{k=1-n}^{n-1}\!\!c_k(f)e^{ikx}$ 
v\'erifie~: $\displaystyle4\pi S^{\ast}_n(f)(x)\!=\!%
\!\int_0^{2\pi}\!\!\!f(t)e^{-in(x-t)}dt+
\int_0^{2\pi}\!\!\!f(t)e^{in(x-t)}dt
+2\!\!\!\sum_{k=1-n}^{n-1}\!\!\int_0^{2\pi}\!\!\!f(t)e^{ik(x-t)}dt$
qui, puisque $e^{ik(x-t)}+e^{-ik(x-t)}\!=\!2\cos(k(x-t))\!=\!2\cos(k(t-x))$,
vaut~:
$\displaystyle\int_0^{2\pi}\!\!\!f(t)\bigl[e^{-in(x-t)}+e^{in(x-t)}+%
2\!\sum_{k=1-n}^{n-1}e^{ik(x-t)}\bigr]dt\!=\!%
\int_0^{2\pi}\!\!\!f(t)2D^{\ast}_n(t-x)dt$. 
D'o\`u le r\'esultat en divisant par $4\pi$ et utilisant la
$2\pi$-p\'eriodicit\'e de
$t\mapsto f(t)D^{\ast}_n(t-x)$.%
\hfill\findem

\nsecp La derni\`ere int\'egration par partie.| 
Si $f$ est
$2\pi$ p\'eriodique, d\'erivable par morceaux avec
$f'$ int\'egrable sur $[0, 2\pi]$ et $x\!\in\!{\bf R\/}$,
on note
$x\!=\alpha_0\!<\!\cdots\!<\!\alpha_{m+1}\!=\!x+2\pi$
tels que $f$ est d\'erivable sur $]\alpha_{j}, \alpha_{j+1}[$.

Ainsi  $f$ a des limites \`a gauche et \`a droite en
les $\alpha_k$ et une int\'egration par parties
de l'int\'egrale de Dirichlet
sur chaque intervalle de d\'erivabilit\'e de $f$
donne~:
$$\displaystyle 2\pi S^{\ast}_n(f)(x)\!=\!\sum_{j=0}^{m}\bigl[%
f_{-}(\alpha_{j+1})E^{\ast}_n(\alpha_{j+1}-x)-%
f_{+}(\alpha_j)E^{\ast}_n(\alpha_j-x)
-\!\int_{\alpha_{j}}^{\alpha_{j+1}}\!\!\!\!\!\!f'(t)E^{\ast}_n(t-x)dt\bigr]$$
qui, en regroupant les termes et comme par p\'eriodicit\'e 
$f_{-}(x+2\pi)=f_{-}(x)$, vaut~:
$$f_{-}(x)E^{\ast}_n(2\pi)-f_{+}(x)E^{\ast}_n(0)%
-\Bigl[\sum_{j=1}^{m}[f_{+}(\alpha_j)-f_{-}(\alpha_j)]E^{\ast}_n(\alpha_j-x)+%
\int_x^{x+2\pi}\!\!\!\!\!\!\!\!\!\!f'(t)E^{\ast}_n(t-x)dt\Bigr]$$
d'o\`u, puisque (Cf. {\bf 3\/}) $E^{\ast}_n(2\pi)\!=\!\pi\!=\!-E^{\ast}_n(0)$, 
et par B
du r\'esum\'e (Cf. {\bf 2\/} et {\bf 3\/})~:
$$\displaystyle
\build{\rm Lim}_{n\rightarrow\infty}^{}S^{\ast}_n(f)(x)\!=\!%
{{f_{-}(x)+f_{+}(x)}\over{2}}$$
\hfill\findem

%% file: Mat.tex
\null
\vskip7mm
\centerline{\petcap  Table des mati\`eres}

\vskip10mm
\centerline{\petcap {\bf I\/} A B C des s\'eries de Fourier%
 }

{ R\'esum\'e.\/}
 \tpage{1}

{  Abstract.\/}
 \tpage {1}

{  Commentaires bibliographiques.\/}
 \tpage {2}

{  R\'ef\'erences.\/}
 \tpage {3}

\vskip10mm
\centerline{\petcap {\bf II\/} Appendices\note{\eightpoint
 $\sum\!\!{{\sin(kx)}\over{k}}^{(\ast)}$\'ecrivait ABKSA :
\`a ne brouter qu'en cas de
s\'echeresse sur l'alpage ({\it N.d.T.\/}$^{(\ast\ast)}$).\hfill\break
$(\ast)$ Chalet de l'Alpe/Habert de Saint Vincent, Longitude 5.91437 Lattitude 45.41463.\hfill\break
$(\ast \ast)$ alexis.charles.marin@gmail.com}%
}%
 
\vskip1mm
{\bf 0\/}   Sommes partielles et sommes partielles modifi\'ees.
 \tpage {4}

{\bf 1\/}   L'expression $D^{\ast}_n(t)\!= \!\sin(nt)\cotg({{t}\over{2}})$
du noyau de Dirichlet modifi\'e.
 \tpage {4}

{\bf  2\/}  La premi\`ere int\'egration par parties.
 \tpage {4}

{\bf  3\/}  La majoration uniforme
$\big{|}E^{\ast}_n(x)\big{|}\!=\!\big{|}\int_{\pi}^xD^{\ast}_n(t)dt\big{|}\leq\pi$.
 \tpage {4}

{\bf 4\/} Fonctions cr\'eneaux en deux coups de la corne $\psi$.
 \tpage {5}

{\bf 5\/} L'int\'egrale de Dirichlet 
$\displaystyle S^{\ast}_n(f)(x)\!=\!%
{{1}\over{2\pi}}\int_x^{x+2\pi}f(t)D^{\ast}_n(t-x)dt$.
 \tpage {5}

{\bf 6\/} La derni\`ere int\'egration par parties.
 \tpage {5}
\vskip10mm

{\petcap Table des mati\`eres.\/}
 \tpage{6}